\documentclass[10pt,twoside]{article}
\usepackage{graphicx}
\usepackage{amsmath}
\usepackage{Latex-document}

\title{\bf  Positively Curved Surfaces \vskip -2mm
in the Three-sphere\vskip 6mm}

\author{B. Andrews\vspace*{-0.5cm}\thanks{Centre for Mathematics and
its Applications, Australian National University, ACT 0200,
Australia. E-mail: andrews@maths.anu.edu.au}}
\date{\vspace{-8mm}}

\begin{document}
\maketitle

\thispagestyle{first} \setcounter{page}{221}

\begin{abstract}\vskip 3mm
In this talk I will discuss an example of the use of fully
nonlinear parabolic flows to prove geometric results.  I will
emphasise the fact that there is a wide variety of geometric
parabolic equations to choose from, and to get the best results it
can be very important to choose the best flow.  I will illustrate
this in the setting of surfaces in a three-dimensional sphere.

There are quite a few relevant results for surfaces in the sphere
satisfying various kinds of curvature equations, including totally
umbillic surfaces, minimal surfaces and constant mean curvature
surfaces, and intrinsically flat surfaces.  Parabolic flows can
strengthen such results by allowing classes of surfaces satisfying
curvature inequalities rather than equalities:  This was first
done by Huisken, who used mean curvature flow to deform certain
classes of surfaces to totally umbillic surfaces.  This motivates
the question ``What is the optimal result of this kind?'' --- that
is, what is the weakest pointwise curvature condition which
defines a class of surfaces which retracts to the space of great
spheres?

The answer to this question can be guessed in view of the
examples.  To prove it requires a surprising choice of evolution
equation, forced by the requirement that the pointwise curvature
condition be preserved.

I will conclude by mentioning some other geometric situations in
which strong results can be proved by choosing the best possible
evolution equation.

\vskip 4.5mm

\noindent {\bf 2000 Mathematics Subject Classification:} 53C44,
53C40.

\noindent {\bf Keywords and Phrases:} Surfaces, Curvature,
Parabolic equations.
\end{abstract}

\vskip 12mm

\section{Introduction} \label{section 1}\setzero
\vskip-5mm \hspace{5mm}

My aim in this talk is to demonstrate the use of fully nonlinear
parabolic evolution equations as tools for proving results in
differential geometry.  I will emphasise the fact that there is a
wide variety of flows which are geometrically defined and
potentially applicable to geometric problems, and that there is
great benefit to be had by choosing the flow carefully.  I will
focus on a  particular application, relating to surfaces in the
3-sphere, but the method has much wider applicability.

There are some well-known examples of geometric evolution
equations of the kind I want to consider:   Eells and Sampson
\cite{ES} used a heat flow to prove existence of harmonic  maps
into non-positively curved targets; Hamilton considered the flow
of Riemannian metrics in the direction of their Ricci tensor, and
proved that it deforms metrics of positive Ricci curvature on
three-manifolds \cite{Ha1} and metrics of positive curvature
operator on four-manifolds \cite{Ha2} to constant curvature
metrics.  The Ricci flow also gives results in higher dimensions,
proved by Huisken \cite{Hu1}, Nishikawa \cite{Ni} and Margerin
\cite{Ma1}--\cite{Ma3}, if the curvature tensor is suitably
pinched.   The mean curvature flow of submanifolds of Euclidean
space is also well-known as the gradient descent flow of the area
functional, and because it arises in models of interfaces such as
in annealing metals. The examples I will concentrate on are
closest to the last example, as they are evolution equations
describing submanifolds moving with curvature-dependent velocity.
There are many parabolic flows of this kind, particularly for the
codimension one (hypersurface) case:  William Firey \cite{F}
introduced the motion by Gauss curvature as a model for pebbles
wearing away as they tumble, and other flows which have been
considered include motion by powers of Gauss curvature \cite{T},
\cite{Ch1}, the square root of the scalar curvature \cite{Ch2},
the harmonic mean of the principal curvatures
\cite{A1}--\cite{A2}, and the reciprocal of the mean curvature
\cite{HI}.  More generally, one can take the velocity to be a
function of the principal curvatures which is monotone increasing
in each argument.

This gives a huge variety of flows to choose from, so it makes
sense to choose the flow carefully to suit the problem.  I will
illustrate a strategy for choosing the flow by asking that some
desired curvature inequality be preserved under the flow.

I will begin, in the next two sections, by discussing some old
results concerning surfaces in the three-sphere.  This motivates
the results of the later sections.

\section{Constant mean curvature surfaces} \label{section 2}
\setzero\vskip-5mm \hspace{5mm}

There is a well-known result of Simons \cite{Si} which says that a
minimal hypersurface in a $S^{n+1}$ with the squared norm of the
second fundamental form $|A|^2$ less than $n$ is in fact totally
geodesic (hence a great $n$-sphere). This result comes from an
application of Simons' identity which relates the second
derivatives of mean curvature to the Laplacian of the second
fundamental form:
\begin{equation*}
\nabla_i\nabla_j H = \Delta h_{ij}
+|A|^2h_{ij}-Hh_i^ph_{pj}+Hg_{ij}-nh_{ij}.
\end{equation*}
  From this we can deduce if the hypersurface is minimal (so $H=0$)
\[
0 = \Delta|A|^2-2|\nabla A|^2 +2|A|^2(|A|^2-n).
\]
If $|A|^2<n$ at a maximum, then the maximum principle implies
$|A|^2$ is identically zero, and the result follows.  Also, if the
maximum of $|A|^2$ is equal to $n$, then $M$ must be a product
$S^k(a)\times S^{n-k}(b)$ in $R^{k+1}\times R^{n+1-k}$, with radii
$a$ and $b$ determined by the fact that $M$ lies in
$S^{n+1}\subset R^{n+2}$ and is minimal.

Simons' argument was taken up by other authors (\cite{Ok},
\cite{CN}, \cite{AdC}) in the slightly more general setting of
constant mean curvature hypersurfaces. The results are similar:
If the hypersurface has constant mean curvature $H$, and $|A|$ is
bounded by a constant depending on $n$ and $H$, then the
hypersurface is totally umbillic, hence a geodesic sphere in
$S^{n+1}$; if the inequality is not strict then the only extra
possibilities are products of spheres. The argument is similar to
that above, but complicated by the non-vanishing of the mean
curvature.

Let me look closer at the situation for surfaces in the
three-sphere:  The intrinsic curvature of the surface is given by
$1+\kappa_1\kappa_2 = 1+\frac12H^2-\frac12|A|^2$.  If $M$ is
minimal, then $H=0$, so $|A|^2<2$ is equivalent to positivity of
the intrinsic curvature.  This is also true for constant mean
curvature surfaces:  In two dimensions, the curvature condition
from  \cite{Ok} and \cite{CN} is equivalent to positivity of the
intrinsic curvature.

%% do Carmo, at al

\section{Flat tori}\label{section 3}
\setzero\vskip-5mm \hspace{5mm}

The condition of positive intrinsic curvature seems natural in
view of the results on constant mean curvature surfaces.  For
surfaces in space, positive curvature is a rather restrictive
condition --- a compact surface satisfying this condition is the
boundary of a convex region.  In the 3-sphere it seems somewhat
less restrictive, as we can see by considering the `boundary' case
of flat surfaces, where there are the beautiful results of Weiner
\cite{W} and Enomoto \cite{E} which classify flat tori in the
3-sphere by their Gauss maps.  It was known for some time that
there are many examples of these (see \cite{P}), since the inverse
image of any smooth curve in $S^2$ under the Hopf projection is a
flat torus in $S^3$.  These examples are all invariant under the
action of $U(1)$ on $C^2\simeq R^4$, but Weiner and Enomoto showed
that there are many examples which are not symmetric.

The Gauss map of a surface in $S^3$ can be thought of in several
ways:  One can consider the tangent plane of the surface as a
subspace of $R^4$, which gives a map from the surface to the
Grassmannian $G_{2,4}$ of 2-planes in $R^4$.  The latter is a
metric product $S^2\times S^2$, and the projections onto each
factor are called the self-dual and anti-self-dual Gauss maps.
Alternatively, since $S^3$ is a group, one can map the unit normal
of the surface by either left or right translations to the Lie
algebra --- this again gives two maps to $S^2$, and of course
these are the same as before:  The self-dual Gauss map is the same
as the left-translation Gauss map, and the anti-self-dual Gauss
map is the same as the right-translation Gauss map.

Enomoto \cite{E} observed that if $M^2$ is intrinsically flat in
$S^3$, then both Gauss maps are degenerate (their images are just
curves in $S^2$).  Weiner gave the complete classification result:
The image curves $\gamma_1$ and $\gamma_2$ necessarily have zero
total curvature, and if $I_1$ and $I_2$ are subintervals of
$\gamma_1$ and $\gamma_2$ respectively, then
$|\int_{I_1}\kappa\,ds|+|\int_{I_2}\kappa\,ds|<\pi$.  Conversely,
if $\gamma_1$ and $\gamma_2$ are any curves satisfying these
conditions, then there is a flat torus with these curves as the
images of the two Gauss maps, and the torus is unique up to motion
by unit speed in the normal direction.

This gives a very large family of flat tori in the 3-sphere, and
from these we see that surfaces with positive intrinsic curvature
in $S^3$ can look quite complicated:  The surface can look
metrically like a long thin cylinder with caps on the ends, placed
in $S^3$ by `winding around' a flat torus many times before
closing off the ends.

\section{Curvature flow} \label{section 4}
\setzero\vskip-5mm \hspace{5mm}

Curvature flow can give powerful generalisations of results like
those from \cite{Si}, \cite{Ok} and \cite{CN}:  Huisken \cite{Hu3}
extended techniques developed earlier for convex hypersurfaces in
Euclidean space \cite{Hu1} to prove the following result:

\smallskip\noindent{\bf Theorem: }\it
Let $M^n_0=x_0(M)$ be a hypersurface in $S^{n+1}$ which satisfies
\[
|A|^2<\frac{1}{n-1}H^2+2
\]
if $n>2$, and
\[
|A|^2<\frac34H^2+\frac43
\]
if $n=2$.  Then there exists a smooth family of hypersurfaces
$\{M_t = x_t(M)\}_{0\leq t<T}$ which satisfy the same curvature
condition and move by mean curvature flow with initial data $M_0$.
Either $T<\infty$ and $M_t$ is asymptotic to a family of geodesic
spheres shrinking to their common centre, or $T=\infty$ and $M_t$
approaches a great sphere. \rm

\smallskip
This includes the result that there are no minimal surfaces with
$|A|^2<n$ except great spheres.  It also implies the stronger
statement that every hypersurface satisfying
$|A|^2<\frac{1}{n-1}H^2+2$ can be deformed, keeping this
condition, to a great sphere (except in the case $n=2$). The
condition $|A|^2<\frac{1}{n-1}H^2+2$ is the same as that arrived
at by Okumura \cite{Ok} for constant mean curvature surfaces
(Cheng and Nakagawa \cite{CN} improved this for higher dimensions,
but in two dimensions it is sharp).  The proof of the above result
is significantly more difficult than that for the constant mean
curvature case.

The result seems very satisfying, except when $n=2$ where the
method does not seem to work for Okumura's condition
$|A|^2<H^2+2$.  The latter is exactly the condition of positive
intrinsic curvature.  This raises several questions:  Does mean
curvature flow in fact preserve this condition?  If not, is there
any flow which does?

\section{The optimal result} \label{section 5}
\setzero\vskip-5mm \hspace{5mm}

\subsection{Choosing the evolution equation}\vskip-5mm \hspace{5mm}

Now we can illustrate the method:  The previous questions can be
answered in a rather systematic way.  The idea is to write down
the conditions required for an arbitrary flow by a function $F$ of
curvature to preserve positive intrinsic curvature.

We can write down an evolution equation for an arbitrary function
$G$ of the principal curvatures $\kappa_1$ and $\kappa_2$, and see
what conditions are required for the flow to preserve the
condition $G\geq 0$. For convenience we can write $G$ in the form
\begin{equation}\label{eq:formG}
G(\kappa_1,\kappa_2) =
(\kappa_1-\kappa_2)^2-\varphi(\kappa_1+\kappa_2)^2
\end{equation}
so that in the case we are interested in,
$\varphi(x)=\sqrt{4+x^2}$. We can also write
\begin{equation}
\label{eq:formF} F=f(\kappa_1+\kappa_2,G).
\end{equation}
Then the evolution equation for $G$ is as follows:
\begin{equation}\label{eq:evolG}
\frac{\partial G}{\partial t} ={\dot F}^{ij}\nabla_i\nabla_jG +
Q(h)(\nabla h,\nabla h) + Z(h),
\end{equation}
where $\dot F$ is the matrix of derivatives of $F$ with respect to
the components of the second fundamental form, which is positive
definite as long as $F$ is an increasing function of each of the
principal curvatures.  The second term is a quadratic function of
the components of the derivative of the second fundamental form,
with coefficients depending on curvature $h$, explicitly given by
\[
Q = \left(\dot G^{ij}\ddot F^{kl,mn}-\dot F^{ij}\ddot
G^{kl,mn}\right) \nabla_ih_{kl}\nabla_jh_{mn},
\]
where $\ddot F$ is the second derivative of $F$ with respect to
the components of $h$. The last term $Z$ depends on the curvature
alone, and has the form
\begin{align*}
Z &= \dot G^{ij}\left(F(h^2_{ij}+g_{ij})+\dot F^{kl}
\left(h_{ij}h^2_{kl}-h_{kl}h^2_{ij}+g_{ij}h_{kl}-g_{kl}h_{ij}
\right)\right)\\
&= F\left(\dot G^1(1+\kappa_1^2)+\dot
G^2(1+\kappa_2^2)\right)+(1+\kappa_1\kappa_2)(\kappa_2-\kappa_1)
(\dot G^1\dot F^2-\dot F^1\dot G^2).
\end{align*}
To show that $G\geq 0$ is preserved (with $G=1+\kappa_1\kappa_2$),
we consider the situation at a point where $G$ first attains a
zero minimum.  Then the first term on the right-hand side of
\eqref{eq:evolG} is non-negative; we consider each of the other
terms.  The last term is simplest: Substituting the forms of $F$
and $G$ from \eqref{eq:formG} and \eqref{eq:formF}, we find
\[
Z = G\left(fH+\frac{\partial f}{\partial H}\varphi^2\right),
\]
so $Z$ vanishes at a zero of $G$, no matter what speed $F$ we use.
This is another indication of the fact that the condition of
positive intrinsic curvature is optimal. The gradient terms are
the most complicated, but we can simplify them significantly by
observing two things:  First, $\nabla h$ is a totally symmetric
$3$-tensor, by the Codazzi equation.  Second, at a minimum of $G$,
the gradients of $G$ vanish.  It follows that there are only two
independent components of $\nabla h$, and one finds that these
never mix in the expression for $Q$, so that
\[
Q = \alpha (\nabla_1h_{22})^2 + \beta (\nabla_2h_{11})^2.
\]
Since we have no further information about $\nabla h$ (that is, no
reason to expect that the magnitudes of these remaining components
should vanish) we must impose the condition that $\alpha$ and
$\beta$ are non-negative.  This gives two conditions, which we can
interpret as conditions on the first and second derivatives of
$F$.  A fact which is perhaps not obvious is that these conditions
only involve the restriction of $F$ to the boundary of the set
$\{G=0\}$ in the curvature plane, so we can consider $F$ as
defined by \eqref{eq:formF} with $G=0$. Then the conditions can be
written explicitly as follows:
\[
\frac{\varphi''}{1-\varphi'}-\frac{1-\varphi'}{\varphi}
\leq\frac{f''}{f'}\leq \frac{\varphi''}{1+\varphi'}
+\frac{1+\varphi'}{\varphi}.
\]
In the case of interest, we have $\varphi = \sqrt{4+H^2}$, and the
first and last quantities are both equal to $-2H/(4+H^2)$.  The
only possibilities for $F$ are the following:
\[
F = C_1 + C_2\arctan\left(\frac{H}2\right).
\]
This applies only along the curve $\{G=0\}$, so we are reasonably
free to choose $F$ in the region where $G>0$, as long as it is
monotone in both principal curvatures.

\subsection{The extreme case}\vskip-5mm \hspace{5mm}

The remarkably restricted form of the evolution equation is
illuminated somewhat by considering the extreme case of flat
surfaces:  If the flow preserves positive intrinsic curvature,
then it must also preserve zero curvature.  As outlined above, the
structure of surfaces with zero curvature is very well understood,
and in particular the Gauss map $G: M^2\to S^2\times S^2$ has the
remarkable property that the projection onto each factor is
one-dimensional.  This must be preserved under the flow.

The flow we have ended up with is characterised by the fact that
the Gauss map evolves according to the mean curvature flow (now
for codimension 2 surfaces in $S^2\times S^2$, which means that
each of the two curves coming from the two projections of the
Gauss map evolves according to the curve-shortening flow in $S^2$.
Since each of the curves divides the area of the sphere into two
equal parts, the image of the Gauss map never develops
singularities (at least in the case where the two curves are
homotopic to great circles traversed once), but in fact the flat
tori will in general develop singularities --- this is analogous
to the motion of a curve in the plane with constant normal speed,
which develops singularities even though the normal direction
stays constant at each point.  Incidentally, there has been some
very impressive recent progress on mean curvature flow in higher
codimension, due to Mu-Tao Wang \cite{Wa1}--\cite{Wa3}, who has
used it to prove several very interesting results regarding maps
between manifolds.

The examples of flat tori can be used to prove that there is no
other curvature-driven flow of surfaces which preserves the
condition of positive curvature, by giving examples for any other
flow of flat tori which do not stay flat.

\subsection{Regularity}\vskip-5mm \hspace{5mm}

% Non-concavity

A technical issue which arises is the following:  The speed we
ended up with is not concave or convex as a function of the second
fundamental form.  The regularity estimates due to Krylov
\cite{Kr} and Evans \cite{Ev} for fully nonlinear equations
(needed to prove that we get classical solutions of the flow)
require concavity, so we cannot use these. Instead it is possible
to adapt the estimates for elliptic equations in two variables
(due to Morrey \cite{Mo} and Nirenberg \cite{Nr}) to give good
$C^{2,\alpha}$ estimates for solutions of fully nonlinear
parabolic equations in two space variables.

\subsection{Curvature pinching}\vskip-5mm \hspace{5mm}

Now we come to the problem of choosing a good way to extend the
speed from the boundary $\{G=0\}$ to the interior of the region
$\{G>0\}$. The idea is to do this in such a way that any compact
surface with strictly positive curvature necessarily has very
strongly controlled curvature in the future --- that is, we want
the region $\{G>0\}$ to be exhausted by a nested family of regions
which stay away from the boundary, and only approach infinity near
the `umbillic' line $\kappa_1=\kappa_2$.  This means that any
singularity which occurs will have to be totally umbillic, so
occurs only when the surface shrinks to a point while becoming
spherical in shape.

This can be done in many ways.  One which is relatively simple to
describe, but results in solutions which are only $C^{2,\alpha}$,
is as follows:  Take
\[
F =\begin{cases}
\arctan\kappa_1+\arctan\kappa_2,&\kappa_1\kappa_2<1;\\
\frac\pi 4(\kappa_1\kappa_2+1),&\kappa_1\kappa_2>1.
\end{cases}
\]
This is then a Lipschitz, monotone increasing function of the
curvatures, and one can check that the following regions of the
curvature plane are preserved:
\[
\Omega_\varepsilon = \left\{|\kappa_1-\kappa_2|\leq\frac{1+
\kappa_1\kappa_2}{\varepsilon}\right\}
\cap\left\{\kappa_1\kappa_2\leq 1\right\}
\cup\left\{|\kappa_1-\kappa_2|\leq\frac2\varepsilon\right\}
\cap\left\{\kappa_1\kappa_2\geq 1\right\}.
\]
This means that the difference between the principal curvatures
stays bounded even if the curvature becomes large, which implies
very strong control on singularities. This is similar to the
estimate used in \cite{A3} to prove that worn stones (i.e. convex
surfaces moving by their Gauss curvature) become round as they
shrink to points.

With a little more work we can choose the speed to be a smooth
function of the principal curvatures, and then solutions are also
smooth.

% Choosing the flow to preserve good inequalities

In the choice above, we also have the nice feature that minimal
surfaces do not move.  We can with slight modifications arrive at
a speed for which constant mean curvature surfaces do not move,
for any particular choice of the mean curvature, as long as we are
willing to work in the category of oriented surfaces.   More
generally, we can contrive that for a given monotone increasing
function $\phi$ of the principal curvatures, surfaces satisfying
$\phi=0$ do not move.   Here $F$ (and $\phi$) must be symmetric.
We can also choose if desired a speed which is always positive, so
that there are no stationary solutions.

\subsection{The results}\vskip-5mm \hspace{5mm}

The main result for the above speed is the following:

\smallskip\noindent{\bf Theorem 1.} \it
Let $x_0$ be an immersion of $S^2$ in $S^3$, with non-negative
intrinsic curvature in the induced metric.  Then the flow
constructed above deforms $M_0=x_0(S^2)$ through a family
$M_t=x_t(S^2)$, with intrinsic curvature strictly positive for
each $t>0$, to either a great sphere (in infinite time) or to a
point, with spherical limiting shape (in finite time).  If $M_0$
is embedded, then so is $M_t$ for each $t>0$. \rm

\smallskip
This includes in particular Simons' result on mimimal surfaces.
If we modify the speed somewhat, then we get the following result,
which gives in particular a new result for Weingarten surfaces in
the 3-sphere:

\smallskip\noindent{\bf Theorem 2. }\it
Let $\phi$ be any smooth, strictly monotone function of $\kappa_1$
and $\kappa_2$ defined on $\{\kappa_1\kappa_2+1\geq 0\}$.  Then
there exists a function $F$ which is smoothly defined on
$\{\kappa_1\kappa_2+1\geq 0\}$, and strictly monotone increasing
in each argument, with $\text{\rm sgn} F=\text{\rm sgn}\phi$
everywhere, such that the following holds:   If $M_0=x_0(S^2)$ is
a smooth compact surface in $S^3$ with non-negative intrinsic
curvature, then the motion with speed $F$ deforms $M_0$ through a
smooth family $\{M_t\}_{0\leq t<T}$, each strictly positively
curved, which either converge to a point with spherical limiting
shape with $T<0$, or converge to a totally umbillic surface
(spherical cap) with $\phi=0$ if $T=\infty$. \rm

\smallskip This includes two cases:  Either there is some point where
$\phi=0$, in which case there is a spherical cap with $\phi=0$ and
the above result implies that this is the only surface with
$\phi=0$ with positive intrinsic curvature, or $\phi$ is never
zero, in which case all surfaces converge to points.  In the
latter case a very small geodesic sphere with one choice of
orientation will shrink inwards to its centre, while the same
sphere with the opposite orientation expands over the equator and
eventually contracts to the antipodal point.  In this way we have
a unique way of associating an oriented surface with the point it
eventually contracts to, and we deduce the following:

\smallskip\noindent{\bf Theorem 3. }\it
The space of oriented surfaces with positive intrinsic curvature
in $S^3$ retracts onto $S^3$. \rm

\smallskip Finally, if we introduce some non-local terms in the speed, we
can devise a flow which fixes the enclosed volume, preserves
positive intrinsic curvature, and gives convergence to spherical
caps, without moving constant mean curvature surfaces.

\section{Other results by related methods} \label{section 6}

\vskip-5mm \hspace{5mm}

The methods I outlined above also yield interesting results for a
variety of other problems:  One which works out similarly, and
which has some interesting parallels, is that of surfaces in
three-dimensional hyperbolic space. The surfaces of interest are
those for which all of the principal curvatures are less than $1$
in magnitude.  We can find a flow which deforms any such surface
in a compact hyperbolic manifold to a minimal surface, while
keeping the principal curvatures less than $1$ in magnitude.
Rather surprisingly, this flow is in a way the hyperbolic analogue
of the one we just described for the sphere:  Instead of moving
with speed equal to the sum of the arctangents of the principal
curvatures, we move with speed equal to the sum of the hyperbolic
arctangents of the principal curvatures.  The resulting flow is
very well-behaved, and has the interesting property that the Gauss
map of the surface (the map which takes a point of the surface to
its tangent plane, thought of as a point in the Grassmannian of
spacelike 2-planes in Minkowksi space $R^{3,1}$), evolves
according to mean curvature flow.

The methods also give good results for hypersurfaces in
higher-dimensional spheres:  Hypersurfaces with positive sectional
curvatures can be deformed in such a way as to preserve that
condition, and similar results can be deduced.  The condition of
positive sectional curvature can probably be relaxed:  Positive
sectional curvature is implied by the condition of Okumura
\cite{Ok} for constant mean curvature hypersurfaces, but not by
the sharper condition of Cheng and Nakagawa \cite{CN} and Alencar
and do Carmo \cite{AdC}.

\label{lastpage}


\begin{thebibliography}{aa}

\bibitem{AdC} H. Alencar and M. do Carmo, Hypersurfaces with constant
mean curvature in spheres, {\it Proc. Amer. Math. Soc.} 120
(1994), 1223--1229.

\bibitem{A1} B. Andrews, Contraction of convex hypersurfaces in
Euclidean space, {\it Calc. Var. P.D.E.} 2 (1994), 151--171.

\bibitem{A2} B. Andrews, Contraction of convex hypersurfaces in
Riemannian spaces, {\it J. Differential Geometry} 39 (1994),
407--431.

\bibitem{A3} B. Andrews, Gauss Curvature Flow: The Fate of the Rolling
Stones, {\it Invent. Math.} 138 (1999), 151--161.

\bibitem{CN} Q.-M. Cheng and H. Nakagawa, Totally umbillic hypersurfaces,
{\it Hiroshima Math. J.} 20 (1990), 1--10.

\bibitem{Ch1} B. Chow, Deforming convex hypersurfaces by the $n$th root
of the Gaussian curvature, {\it J. Differential Geom.} 22 (1985),
117--138.

\bibitem{Ch2} B. Chow, Deforming convex hypersurfaces by the square root
of the scalar curvature, {\it Invent. Math.} 87 (1987), 63--82.

\bibitem{ES} J. Eells and J. Sampson, Harmonic mappings of Riemannian
manifolds, {\it Amer. J. Math.} 86 (1964), 109--160.

\bibitem{E} K. Enomoto,
The Gauss image of flat surfaces in $R\sp 4$, {\it Kodai Math. J.}
9 (1986), 19--32.

\bibitem{Ev} L. C. Evans, Classical solutions of fully
nonlinear, convex, second order elliptic equations, {\it
Comm.~Pure Appl.~Math.} 24 (1982), 333--363.

\bibitem{F} W. J. Firey, Shapes of worn stones.
{\it Mathematika} 21 (1974), 1--11.

\bibitem{Ha1} R. S. Hamilton, Three-manifolds with positive Ricci
curvature, {\it J. Differential Geometry}, 17 (1982), 255--306.

\bibitem{Ha2} Four-manifolds with positive curvature operator, {\it J.
Differential Geometry} 24 (1986), 153--179.

\bibitem{Hu1} G. Huisken, Flow by mean curvature of convex surfaces into
spheres, {\it J. Differential Geometry} 20 (1984), 237--266.

\bibitem{Hu2} G. Huisken, Ricci deformation of the metric on a
Riemannian manifold, {\it J. Differential Geometry} 21 (1985),
47--62.

\bibitem{Hu3} G. Huisken, Deforming hypersurfaces of the sphere by their
mean curvature, {\it Math. Z.} 195 (1987), 205--219.

\bibitem{HI} G. Huisken and T. Ilmanen, The Riemannian Penrose
Inequality, {\it Internat. Math. Res. Notices} 1997, no. 20,
1045--1058.

\bibitem{Kr} N. V. Krylov, Boundedly inhomogeneous elliptic and
parabolic equations, {\it  Izvestia Akad. Nauk. SSSR} 46 (1982),
487--523.  English translation in {\it Math. USSR Izv.}  20
(1983).

\bibitem{Ma1} C. Margerin, Pointwise pinched manifolds are space forms,
{\it Proc. Symp. Pure Math.} 44, 1986.

\bibitem{Ma2} C. Margerin, Une caractrisation optimale de la structure
diffrentielle standard de la sphre en terme de courbure pour
(presque) toutes les dimensions, {\it C. R. Acad. Sci. Paris S\'er
I Math.} 319 (1994) 713--716 and 605--607.

\bibitem{Ma3} C. Margerin, A sharp characterization of the smooth
$4$-sphere in curvature terms, {\it Comm. Anal. Geom.} 6 (1998),
21--65.

\bibitem{Mo} C.B.~Morrey, Jr.,  On the solutions of quasi-linear
elliptic partial differential equations, {\it
Trans.~Amer.~Math.~Soc.} 43, (1938), 126--166.

\bibitem{Nr} L. Nirenberg, On nonlinear elliptic partial
differential equations and H\"older continuity, {\it Comm.~Pure
Appl.~Math.}  6 (1953),  103--156.

\bibitem{Ni} S. Nishikawa, Deformation of Riemannian metrics and
manifolds with bounded curvature ratios, {\it Proc. Sympos. Pure
Math.} 44, 1986.

\bibitem{Ok} M. Okumura, Hypersurfaces and a pinching problem on
the second fundamental tensor, {\it Amer. J. Math.} 96 (1974),
207--213.

\bibitem{P} U. Pinkall, Hopf Tori in $S^3$, {\it Invent. Math.} 81
(1985), 379--386.

\bibitem{Si} J. Simons, Minimal varieties in Riemannian manifolds,
{\it Ann. of Math.} (2) 88 (1968), 62--105.

\bibitem{T} Kaising Tso, Deforming a hypersurface by its Gauss-Kronecker
curvature, {\it Comm. Pure Appl. Math.} 38 (1985), 867--882.

\bibitem{Wa1} M.-T. Wang, Mean curvature flow of surfaces in Einstein
four-manifolds, {\it J. Differential Geom.} 57 (2001), 301--338.

\bibitem{Wa2} M.-T. Wang, Deforming area preserving diffeomorphism of
surfaces by mean curvature flow, {\it Math. Res. Lett.} 8 (2001),
651--661.

\bibitem{Wa3} M.-T. Wang, Subsets of Grassmannians preserved by mean
curvature flow, preprint, 2002.

\bibitem{W} J. Weiner, Flat tori in $S^3$ and their Gauss maps, {\it
Proc. London Math. Soc.} 62 (1991), 54--76.


\end{thebibliography}
\end{document}